\documentclass[11pt,reqno]{amsart}
\usepackage{amssymb,latexsym}
\usepackage{graphicx}
\usepackage{url}		
\usepackage{tikz}
\usepackage{float}

\calclayout

\makeatletter
\newcommand{\monthyear}[1]{%
  \def\@monthyear{\uppercase{#1}}}
\newcommand{\volnumber}[1]{%
  \def\@volnumber{\uppercase{#1}}}
\AtBeginDocument{%
\def\ps@plain{\ps@empty
  \def\@oddfoot{\@monthyear \hfil \thepage}%
  \def\@evenfoot{\thepage \hfil \@volnumber}}
\def\ps@firstpage{\ps@plain}
\def\ps@headings{\ps@empty
  \def\@evenhead{%
    \setTrue{runhead}%
    \def\thanks{\protect\thanks@warning}%
    \uppercase{}\hfil}%
  \def\@oddhead{%
    \setTrue{runhead}%
    \def\thanks{\protect\thanks@warning}%
    \hfill\uppercase{Visualize Geometric Series}}%
  \let\@mkboth\markboth
  \def\@evenfoot{%
    \thepage \hfil \@volnumber}%
  \def\@oddfoot{%
    \@monthyear \hfil \thepage}%
  }%
\footskip=25pt
\pagestyle{headings}%
}
\makeatother

\theoremstyle{plain}
\numberwithin{equation}{section}
\newtheorem{thm}{Theorem}[section]

\newtheorem{lemma}[thm]{Lemma}

\begin{document}
\monthyear{}
\volnumber{}
\setcounter{page}{1}

\title{Visualize Geometric Series}

\author{H\`ung Vi\d{\^e}t Chu}
\address{Department of Mathematics\\
                Washington and Lee University\\
                Lexington\\
                VA 24450, USA}
\email{hchu@wlu.edu}

\begin{abstract}
We review Mabry's, Edgar's, and the Viewpoints 2000 Group's proofs without words for the geometric series formula.
    Mabry and Edgar proved without words that
    $$\frac{1}{4} + \left(\frac{1}{4}\right)^2 + \left(\frac{1}{4}\right)^3 + \cdots\ =\ \frac{3}{4}\quad\mbox{ and }\quad\frac{4}{9} + \left(\frac{4}{9}\right)^2 + \left(\frac{4}{9}\right)^3 + \cdots\ =\ \frac{4}{5},$$
    respectively. 
    We show that their proofs satisfy certain requirements that make them unique. We then illustrate a common idea between their and the Viewpoints 2000 Group's proofs.
\end{abstract}

\maketitle

\section{Introduction}
Probably all calculus students have encountered the following formula to compute the sum of the geometric series $a + ar + ar^2 + ar^3+ \cdots $ when $|r| < 1$: 
\begin{align}\label{k00}a + ar + ar^2 + ar^3 + \cdots \ =\ \frac{a}{1-r}.\end{align}
It is one of a few precious tools to find the exact sum of a convergent series. Other well-known tools include telescoping series and the Weierstrass factorization theorem, used by Euler in computing $\sum_{n=1}^\infty 1/n^2$. It is worth emphasizing that \eqref{k00} computes the sum of all geometric series as long as the common ratio has a magnitude less than $1$, but  its usefulness does not stop there. For example, differentiating both sides of \eqref{k00} with respect to $r$ gives another formula:
$$\sum_{n=1}^\infty anr^{n-1}\ =\ a + 2ar + 3ar^2 + \cdots\ =\ \frac{a}{(1-r)^2}, \mbox{ for }|r| < 1.$$

While \eqref{k00} is simply the Maclaurin series for $a(1-r)^{-1}$, a more elegant proof is to find a closed form for the partial sums $s_n$ of the geometric series and then let $n$ go to infinity. In particular, observe that for $n\ge 1$, 
\begin{align*}
&(1-r)(a+ar+ar^2+\cdots+ar^{n-1})\\
\ =\ &(a+ar+ar^2+\cdots+ar^{n-1}) - (ar+ar^2+ar^3+\cdots+ar^n)\ =\ a - ar^n.
\end{align*}
Hence,
$$s_n\ :=\ a+ar+ar^2+\cdots+ar^{n-1} \ =\ \frac{a-ar^n}{1-r}.$$
Letting $n\rightarrow \infty$ and noting that $\lim_{n\rightarrow \infty} r^{n} = 0$ because $|r|<1$, we obtain the desired formula. Though the proof is elegant and explanatory, there is still the need for more intuition behind the formula. As a result, many clever proofs without words have been devised. 

In 1999, Mabry \cite{Ma} provided Figure \ref{fig:Mabry's proof} as a visual proof of 
\begin{equation}\label{re1}\frac{1}{4} + \left(\frac{1}{4}\right)^2+ \left(\frac{1}{4}\right)^3+\cdots \  =\ \frac{1}{3}.\end{equation}

\begin{center}
\begin{figure}[H]
\begin{tikzpicture}
[acteur/.style={circle, fill=black,thick, inner sep=0.3pt, minimum size=0.01cm}] 
\draw (-2,0)--(2,0);
\draw (-2,0)--(0,3.4641);
\draw (2,0)--(0,3.4641);
\filldraw[draw=black, fill=yellow] (-1,1.732) -- (1, 1.732) -- (0,0) -- cycle;
\filldraw[draw=black, fill=yellow] (-0.5,2.598) -- (0, 1.732) -- (0.5,2.598) -- cycle;
\filldraw[draw=black, fill=yellow] (0.25,3.031) -- (0, 2.598) -- (-0.25,3.031) -- cycle;
\filldraw[draw=black, fill=yellow] (0.125, 3.2476) -- (0, 3.031) -- (-0.125, 3.2476) -- cycle;
\filldraw[draw=black, fill=yellow] (0.0625, 3.356) -- (0, 3.2476) -- (-0.0625, 3.356) -- cycle;
\filldraw[draw=black, fill=yellow] (0.03125, 3.41) -- (0, 3.356) -- (-0.03125, 3.41) -- cycle;
\draw (2,0)--(1, 1.732) node[pos=0.6,right] {\tiny{1\textsuperscript{st} layer}};
\draw (1,1.732)--(0.5, 2.598) node[pos=0.6, right] {\tiny{2\textsuperscript{nd} layer}};
\draw (0.5,2.598)--(0,3.464) node[pos=0.3, right] {\tiny{3\textsuperscript{rd} layer}};
\end{tikzpicture}
\caption{Mabry's proof.} \label{fig:Mabry's proof}
\end{figure}
\end{center}
Though Mabry's proof is only for the case $a = r = 1/4$ in \eqref{k00}, it provides a nice intuition. If the area of the outermost triangle is $1$, then since the shaded triangles cover a third of each layer, their total area is $1/3$. Another way to compute the total shaded area is by summing up individual triangles: 
$1/4 + (1/4)^2 + (1/4)^3 + \cdots$. The two ways must give the same answer, so we have \eqref{re1}. 

With a similar idea, Edgar \cite{Ed} used Figure \ref{fig: Edgar's} to prove
$$\frac{4}{9}+\left(\frac{4}{9}\right)^2+\left(\frac{4}{9}\right)^3 + \cdots\ =\ \frac{4}{5},$$
which is \eqref{k00} with $a = r = 4/9$.

\begin{center}
\begin{figure}[H]
\begin{tikzpicture}
[acteur/.style={circle, fill=black,thick, inner sep=0.3pt, minimum size=0.01cm}] 
\draw (-2,0)--(2,0);
\draw (-2,0)--(0,3.4641);
\draw (2,0)--(0,3.4641);
\filldraw[draw=black, fill=yellow] (-2,0) -- (-1.333, 1.155) -- (-0.666,0) -- cycle;
\filldraw[draw=black, fill=yellow] (-0.666,0) -- (-1.333, 1.155) -- (0,1.155) -- cycle;
\filldraw[draw=black, fill=yellow] (0.666,0) -- (0, 1.155) -- (1.333,1.155) -- cycle;
\filldraw[draw=black, fill=yellow] (0.666,0) -- (1.333, 1.155) -- (2,0) -- cycle;
\filldraw[draw=black, fill=cyan] (-1.333,1.155) -- (-0.888, 1.925) -- (-0.444,1.155) -- cycle;
\filldraw[draw=black, fill=cyan] (-0.888, 1.925) -- (-0.444, 1.155) --(0, 1.925)--cycle;
\filldraw[draw=black, fill=cyan] (0.888, 1.925) -- (0.444, 1.155) --(0, 1.925)--cycle;
\filldraw[draw=black, fill=cyan] (1.333,1.155) -- (0.888, 1.925) -- (0.444,1.155) -- cycle;
\filldraw[draw=black, fill=green] (-0.888,1.925) -- (-0.593, 2.438) -- (-8/27,1.925) -- cycle;
\filldraw[draw=black, fill=green] (-0.593, 2.438) -- (-8/27,1.925) --(0, 2.438)--cycle;
\filldraw[draw=black, fill=green] (0.593, 2.438) -- (8/27,1.925) --(0, 2.438)--cycle;
\filldraw[draw=black, fill=green] (0.888,1.925) -- (0.593, 2.438) -- (8/27,1.925) -- cycle;
\filldraw[draw=black, fill=red] (-16/27, 2.438) -- (-32/81, 2.78) -- (-16/81, 2.438) -- cycle;
\filldraw[draw=black, fill=red] (-32/81, 2.78) -- (-16/81, 2.438) --(0, 2.78)--cycle;
\filldraw[draw=black, fill=red] (32/81, 2.78) -- (16/81, 2.438) --(0, 2.78)--cycle;
\filldraw[draw=black, fill=red] (16/27, 2.438) -- (32/81, 2.78) -- (16/81, 2.438) -- cycle;
\filldraw[draw=black, fill=orange] (-32/81, 2.78) -- (-64/243, 3) -- (-32/243, 2.78) -- cycle;
\filldraw[draw=black, fill=orange] (-64/243, 3) -- (-32/243, 2.78) -- (0,3) -- cycle;
\filldraw[draw=black, fill=orange] (64/243, 3) -- (32/243, 2.78) -- (0,3) -- cycle;
\filldraw[draw=black, fill=orange] (32/81, 2.78) -- (64/243, 3) -- (32/243, 2.78) -- cycle;
\filldraw[draw=black, fill=gray] (-64/243, 3) -- (-128/729, 3.152) -- (-64/729, 3) -- cycle;
\filldraw[draw=black, fill=gray] (-128/729, 3.152) -- (-64/729, 3) -- (0,3.152) -- cycle;
\filldraw[draw=black, fill=gray] (128/729, 3.152) -- (64/729, 3) -- (0,3.152) -- cycle;
\filldraw[draw=black, fill=gray] (64/243, 3) -- (128/729, 3.152) -- (64/729, 3) -- cycle;
\node (a9) at (0,3.25) [acteur][label = below: ]{};
\node (a10) at (0,3.35) [acteur][label = below: ]{};
\node (a11) at (0,3.3) [acteur][label = below: ]{};
\end{tikzpicture}
\caption{Edgar's proof.} \label{fig: Edgar's}
\end{figure}
\end{center}

There is a small difference between Figures \ref{fig:Mabry's proof} and \ref{fig: Edgar's}. Figure \ref{fig:Mabry's proof}
is divided into four triangles, then we add successive subdivisions to the upper $1/4$. Figure 2 is divided into
nine triangles but we remove existing lines from the upper $4/9$ to divide that into nine triangles, then again leave only
the bottom layer, and repeat the process.

Edgar asked, ``Is it possible to determine which other series allow analogous proof without words?''. The answer is, ``It depends on what we mean by ``analogous" proofs!". If the question asked for a series with $r=a$ as in Mabry's and Edgar's proofs, then we shall see that Figures \ref{fig:Mabry's proof} and \ref{fig: Edgar's} are unique of their kinds; if $r$ is allowed to be different from $a$, then there are many others. 

We conclude this section with the beautiful proof without words by the Viewpoints 2000 Group \cite{VP}: Figure \ref{fig: Viewpoints 2000 Group} proves \eqref{k00} for all $r\in (0,1)$ rather than for special cases as Figures \ref{fig:Mabry's proof} and \ref{fig: Edgar's}, and $a$ and $r$ are not necessarily equal; however, we shall show that all of them share a common proof idea.

\begin{center}
\begin{figure}[H]
\begin{tikzpicture}
[acteur/.style={circle, fill=black,thick, inner sep=0.3pt, minimum size=0.01cm}] 
\draw (0,0)--(5,0) node[pos = 1, below] {$\scriptscriptstyle{x}$};
\draw (0,0)--(0,5) node[pos=0.15, left, xshift=3pt] {$\scriptscriptstyle{a}$};
\draw (0,0)--(0,5) node[pos=1, left] {$\scriptscriptstyle{y}$};
\draw (0,0)--(4.5,4.5) node[pos=0.15, right] {$\scriptscriptstyle{y=x}$};;
\draw (0,1.5)--(4.5,4.5) node[pos = 0.25, left] {$\scriptscriptstyle{y=rx + a}$};
\draw (0,1.5)--(1.5,1.5) node[pos=0.5, below] {$\scriptscriptstyle{a}$};
\draw (1.5, 1.5) -- (1.5, 2.5) node[pos=0.5, left, xshift=3pt] {$\scriptscriptstyle{ar}$}; 
\draw (1.5, 2.5)--(2.5, 2.5) node[pos=0.5, below, yshift=3pt] {$\scriptscriptstyle{ar}$};
\draw (2.5, 2.5)--(2.5, 3.166);
\draw (2.5, 3.166)--(3.166, 3.166);
\draw (3.166, 3.166)--(3.166, 3.61);
\draw (3.166, 3.61)--(3.61, 3.61);
\draw (3.61, 3.61)--(3.61, 3.91);
\draw (3.61, 3.91) -- (3.91, 3.91);
\draw (3.91, 3.91) -- (3.91, 4.11);
\draw (3.91,4.11)--(4.11, 4.11);
\draw (4.11, 4.11)--(4.11, 4.24);
\draw (4.11, 4.24)--(4.24, 4.24);
\draw (4.24, 4.24)--(4.24, 4.33);
\draw (4.24, 4.33)--(4.33, 4.33);
\draw (4.33, 4.33)--(4.33, 4.39);
\draw (4.33, 4.39)--(4.39, 4.39);
\draw (4.5, 4.5)--(4.5,4.5) node[pos = 0.5, right] {$\scriptscriptstyle{\left(\frac{a}{1-r},\frac{a}{1-r}\right)}$};
\end{tikzpicture}
\caption{Viewpoints 2000 Group's proof: 
\begin{equation}\label{re20} a + ar + ar^2 + \cdots \ =\ \frac{a}{1-r},\mbox{ for } 0 < r < 1.\end{equation}} \label{fig: Viewpoints 2000 Group}
\end{figure}
\end{center}
The clever idea of Figure \ref{fig: Viewpoints 2000 Group} lies in the slopes of the two lines $y=x$ and $y = rx+a$ so that moving horizontally preserves the previous distance, while moving vertically shrinks the previous distance by a factor of $r$. Hence, the sum of vertical distances produces the wanted series, $a+ar+ar^2+\cdots$.

\section{Generalizing Mabry's and Edgar's proof}
We shall analyze Mabry's and Edgar's proofs to generalize them to other values of $a$ and $r$. Interestingly, their proofs are the only of their kinds if we ask for the following properties:
\begin{itemize}
\item[(P1)] Parallel line segments with the triangle's base partition the outermost triangle into adjacent layers. Every layer consists of the same number of equal equilateral triangles, the same number of which are shaded in every layer. Furthermore, the height of these triangles is equal to the height of the layer that contains them. 
\item[(P2)] The lengths of the parallel line segments including the triangle's base form a geometric progression. In Mabry's and Edgar's proof, consecutive lengths maintain a constant ratio of $1/2$ and $2/3$, respectively.
\item[(P3)]  The first term of the geometric series to be computed is equal to the common ratio, i.e., $r = a$.
\end{itemize}

Figure \ref{fig: backbones} is the skeleton of Figures \ref{fig:Mabry's proof} and \ref{fig: Edgar's}. The area of $\triangle{\textbf{ABC}} = 1$, thus the length $|\textbf{BC}| = 2\cdot 3^{-1/4}$. (Can you use the geometry of an equilateral triangle to see why this is true?)
Let $n$ and $k$ be the number of triangles and of shaded triangles in each layer, respectively. Then $1\le k<n$. Let $s\in (0,1)$ be the ratio of the lengths of consecutive parallel line segments. For example, $s$ is $1/2$ and $2/3$ in Figures \ref{fig:Mabry's proof} and \ref{fig: Edgar's}, respectively. 
\begin{center}
\begin{figure}[H]
\begin{tikzpicture}
[acteur/.style={circle, fill=black,thick, inner sep=0.3pt, minimum size=0.01cm}] 
\draw (-2,0)--(2,0) node[pos=1, right, xshift=-3pt] {$\scriptscriptstyle{\mathbf B}$};
\draw (-2,0)--(2,0) node[pos=0, left, xshift=3pt] {$\scriptscriptstyle{\mathbf C}$};
\draw (-2,0)--(2,0) node[pos=0.5, above, yshift=-3.5pt] {$\scriptscriptstyle{2\cdot 3^{-1/4}}$};
\draw (-2,0)--(0,3.4641);
\draw (2,0)--(0,3.4641) node[pos=1, above, yshift=-3pt] {$\scriptscriptstyle{\mathbf A}$};
\draw (1.6, 0.693)--(-1.6,0.693) node[pos=0.5, above, yshift=-3.5pt] {$\scriptscriptstyle{2\cdot 3^{-1/4}s}$};
\draw (1.6, 0.693)--(-1.6,0.693) node[pos=0, right, xshift=-2pt] {$\scriptscriptstyle{\mathbf D}$};
\draw (1.6, 0.693)--(-1.6,0.693) node[pos=1, left, xshift=2pt] {$\scriptscriptstyle{\mathbf E}$};
\draw (1.28, 1.247)--(-1.28, 1.247) node[pos=0.5, above, yshift=-3.5pt] {$\scriptscriptstyle{2\cdot 3^{-1/4}s^2}$};
\draw (1.024, 1.69)--(-1.024, 1.69) node[pos=0.5, above, yshift=-3.5pt] {$\scriptscriptstyle{2\cdot 3^{-1/4}s^3}$};
\draw (0.819, 2.045)--(-0.819, 2.045);
\draw (0.655, 2.329)--(-0.655, 2.329);
\draw (0.524,2.556)--(-0.524,2.556);
\draw (1.28, 1.247)--(1.024,1.69) node[pos=0.5, right] {\tiny{3\textsuperscript{rd} layer}};
\draw (1.28, 1.247)--(1.6,0.693) node[pos=0.5, right] {\tiny{2\textsuperscript{nd} layer}};
\draw (1.6, 0.693)--(2,0) node[pos=0.5, right] {\tiny{1\textsuperscript{st} layer}};
\filldraw[draw=black, fill=cyan] (1.2,0) -- (1.6, 0.693) -- (2,0) -- cycle;
\filldraw[draw=black, fill=cyan] (0.96,0.6928) -- (1.28, 1.247) -- (1.6, 0.693) -- cycle;
\filldraw[draw=black, fill=cyan] (0.768, 1.247) -- (1.024, 1.69) -- (1.28, 1.247) -- cycle;
\filldraw[draw=black, fill=cyan] (0.614, 1.69) -- (0.819, 2.045) -- (1.024, 1.69) -- cycle;
\filldraw[draw=black, fill=cyan] (0.491, 2.045) -- (0.655, 2.329) -- (0.819, 2.045) -- cycle;
\filldraw[draw=black, fill=cyan] (0.524,2.556) -- (0.655, 2.329) -- (0.393,2.329) -- cycle;
\node (a9) at (0,2.8) [acteur][label = below: ]{};
\node (a10) at (0,3.2) [acteur][label = below: ]{};
\node (a11) at (0,3) [acteur][label = below: ]{};
\end{tikzpicture}
\caption{An abstract of Mabry's and Edgar's proofs.} \label{fig: backbones}
\end{figure}
\end{center}

We compute the shaded area, denoted by $\mathbf{T}$, in two ways. In each layer, there are $n$ triangles, $k$ of which are shaded, so $\mathbf{T} = k/n$. 

Another way to compute $\mathbf{T}$ is to sum up the areas of the shaded triangles. Since $|\mathbf{DE}| = s|\mathbf{BC}|$, the area of $\triangle \textbf{ADE} = s^2\triangle \textbf{ABC} = s^2$. Hence, the area of the $1$\textsuperscript{st} layer is $1-s^2$, meaning each triangle in the $1$\textsuperscript{st} layer has area 
$u := (1-s^2)/n$. 

Property (P2) guarantees that the edge ratio of two triangles from adjacent layers is $s$. Thus, moving from one layer to the immediate layer above shrinks the triangle area by $s^2$. It follows that each triangle in the  $j$\textsuperscript{th} layer has area $us^{2(j-1)}$. 
Therefore, the total shaded area is (recall that in each layer, $k$ triangles are shaded)
$$ \mathbf{T} \ =\ ku+kus^2+kus^4+\cdots\ =\ ku\sum_{j=0}^\infty s^{2j}.$$
As the two ways of computing $\mathbf T$ must produce the same result, we have the following identity.

\begin{lemma} For integers $1\le k <n$ and $s\in (0,1)$. If $u = (1-s^2)/n$, then
\begin{align}\label{k1}
    ku+kus^2+kus^4+ \cdots \ =\ \frac{k}{n}.
\end{align}
\end{lemma}

We now employ \eqref{k1} to obtain the two identities proven by Mabry and Edgar. In Mabry's proof, we have $(n,k,s) = (3, 1, 1/2)$. Then $u = 1/4$, and \eqref{k1} gives
$$\frac{1}{4} + \left(\frac{1}{4}\right)^2 + \left(\frac{1}{4}\right)^3 + \cdots \ =\ \frac{1}{3}.$$
In Edgar's proof, we have $(n, k, s) = (5,4,2/3)$. Then $u = 1/9$, and \eqref{k1} gives
$$\frac{4}{9} + \left(\frac{4}{9}\right)^2 + \left(\frac{4}{9}\right)^3 + \cdots \ =\ \frac{4}{5}.$$

We are now ready to prove our main result.

\begin{thm}
Mabry's and Edgar's proofs are the only ones that satisfy all the Properties \textnormal{(P1)}, \textnormal{(P2)}, and \textnormal{(P3)}.    
\end{thm}

\begin{proof}
We shall show that Mabry's and Edgar's proofs require two conditions on $s$, which are simultaneously satisfied only when $s$ is $1/2$ or $2/3$.

First, to satisfy Property (P3), we need 
$$kus^2 \ =\ k^2u^2,$$
which implies
\begin{align}\label{k5}\frac{s^2}{1-s^2}\ =\ \frac{k}{n}.\end{align}
Since $k/n < 1$, we know that $0 < s < 1/\sqrt{2}$.
This is the first restriction on $s$.

That $|\textbf{DE}| = s|\textbf{BC}|$ implies that the height of a triangle in the first layer is $(1-s)$ times the height of $\triangle \textbf{ABC}$. Hence, each triangle in the $1$\textsuperscript{st} layer has area $(1-s)^2$, while the $1$\textsuperscript{st} layer has area $1-s^2$. Therefore, Property (P1) asks that $n = (1-s^2)/(1-s)^2\in \mathbb{N}$, so $$\frac{2}{1-s}\ =\ \frac{1-s^2}{(1-s)^2} + 1\in \mathbb{N}.$$ Write $s = 1-2/m$ for some integer $m\ge 3$. Then the number of triangles in each layer is $n = m-1$.
Plugging this back to \eqref{k5} gives
\begin{align*}
    \frac{(1-2/m)^2}{1-(1-2/m)^2} \ =\ \frac{k}{m-1}.
\end{align*}
We obtain $k = m^2/4+1-m$. Hence, $m^2/4$ is an integer, implying that $m = 2m'$ for some  $m'\in \mathbb{N}$. As a result, $s = 1-1/m'$.  This is the second restriction on $s$. 

The two restrictions on $s$ imply that $s\in \{1/2, 2/3\}$. We conclude that Mabry's and Edgar's proofs are the only ones that satisfy all of the Properties (P1), (P2), and (P3). 
\end{proof}

If we drop Property (P3), then we have many other similar proofs to Figures \ref{fig:Mabry's proof} and \ref{fig: Edgar's}. From the above analysis, we need only to choose $s = 1-2/m$ for some $m\ge 3$ and partition each layer into $m-1$ equal triangles. However, this is only true when $m$ is even. Figure \ref{notfit} considers the case $m=3$, a representative for odd $m$.

\begin{center}
\begin{figure}[H]
\begin{tikzpicture}
[acteur/.style={circle, fill=black,thick, inner sep=0.3pt, minimum size=0.01cm}] 
\draw (-2,0)--(2,0);
\draw (-2,0)--(0,3.4641);
\draw (2,0)--(0,3.4641);
\draw (-2/3, 2.309)--(2/3, 2.309);
\draw (2/3, 2.309)--(-2/3,0);
\draw (-2/9,3.079)--(2/9, 3.079);
\draw (2/9, 3.079)--(-2/9, 2.309);
\draw (-2/27, 3.336)--(2/27, 3.336);
\draw (2/27, 3.336)--(-2/27, 3.079);
\end{tikzpicture}
\caption{The case $m=3$. No layer can be partitioned into two equal triangles.} \label{notfit}
\end{figure}
\end{center}

To fix this problem, Figure \ref{fig:m=3} partitions each layer into right (instead of equilateral) triangles.

\begin{center}
\begin{figure}[H]
\begin{tikzpicture}
[acteur/.style={circle, fill=black,thick, inner sep=0.3pt, minimum size=0.01cm}] 
\draw (-2,0)--(2,0);
\draw (-2,0)--(0,3.4641);
\draw (2,0)--(0,3.4641);
\filldraw[draw=black, fill=yellow] (-2,0) -- (-2/3, 2.309) -- (-2/3,0) -- cycle;
\filldraw[draw=black, fill=yellow] (-2/3,0) -- (2/3, 2.309) -- (2/3,0) -- cycle;
\draw (-2/3, 2.309)--(2/3, 2.309);
\filldraw[draw=black, fill=yellow] (-2/3,2.309) -- (-2/9, 3.079) -- (-2/9,2.309) -- cycle;
\filldraw[draw=black, fill=yellow] (-2/9,2.309) -- (2/9, 3.079) -- (2/9,2.309) -- cycle;
\draw (-2/9,3.079)--(2/9, 3.079);
\filldraw[draw=black, fill=yellow] (-2/9,3.079) -- (-2/27, 3.336) -- (-2/27,3.079) -- cycle;
\filldraw[draw=black, fill=yellow] (-2/27,3.079) -- (2/27, 3.336) -- (2/27,3.079) -- cycle;
\draw (-2/27, 3.336)--(2/27, 3.336);
\end{tikzpicture}
\caption{Proof that $\frac{4}{9} + \frac{4}{9}\left(\frac{1}{9}\right) + \frac{4}{9}\left(\frac{1}{9}\right)^2 + \frac{4}{9}\left(\frac{1}{9}\right)^3 + \cdots =  \frac{1}{2}$.} \label{fig:m=3}
\end{figure}
\end{center}

Let us consider $m = 8$. 

\begin{center}
\begin{figure}[H]
\begin{tikzpicture}
[acteur/.style={circle, fill=black,thick, inner sep=0.3pt, minimum size=0.01cm}] 
\draw (-2,0)--(2,0);
\draw (-2,0)--(0,3.4641);
\draw (2,0)--(0,3.4641);
\filldraw[draw=black, fill=yellow] (-2, 0)--(-3/2, 0.866)--(-1, 0)--cycle;
\filldraw[draw=black, fill=yellow] (-1, 0)--(0, 0)--(-1/2, 0.866)--cycle;
\filldraw[draw=black, fill=yellow] (1, 0)--(0, 0)--(1/2, 0.866)--cycle;
\filldraw[draw=black, fill=yellow] (2, 0)--(3/2, 0.866)--(1, 0)--cycle;
\draw (-3/2,0.866)--(3/2, 0.866);
\filldraw[draw=black, fill=yellow] (-3/2, 0.866)--(-9/8, 1.516)--(-3/4, 0.866)--cycle;
\filldraw[draw=black, fill=yellow] (-3/4, 0.866)--(-3/8, 1.516)--(0, 0.866)--cycle;
\filldraw[draw=black, fill=yellow] (3/4, 0.866)--(3/8, 1.516)--(0, 0.866)--cycle;
\filldraw[draw=black, fill=yellow] (3/2, 0.866)--(9/8, 1.516)--(3/4, 0.866)--cycle;
\filldraw[draw=black, fill=yellow] (-9/8, 1.516)--(-27/32, 2)--(-9/16, 1.516)--cycle;
\filldraw[draw=black, fill=yellow] (-9/16, 1.516)--(-9/32, 2)--(0, 1.516)--cycle;
\filldraw[draw=black, fill=yellow] (9/16, 1.516)--(9/32, 2)--(0, 1.516)--cycle;
\filldraw[draw=black, fill=yellow] (9/8, 1.516)--(27/32, 2)--(9/16, 1.516)--cycle;
\filldraw[draw=black, fill=yellow] (-27/32, 2)--(-81/128, 2.366)--(-27/64, 2)--cycle;
\filldraw[draw=black, fill=yellow] (-27/64, 2)--(-27/128, 2.366)--(0, 2)--cycle;
\filldraw[draw=black, fill=yellow] (27/64, 2)--(27/128, 2.366)--(0, 2)--cycle;
\filldraw[draw=black, fill=yellow] (27/32, 2)--(81/128, 2.366)--(27/64, 2)--cycle;
\filldraw[draw=black, fill=yellow] (-81/128, 2.366)--(-0.474, 2.641)--(-0.316, 2.366)--cycle;
\filldraw[draw=black, fill=yellow] (-81/256, 2.366)--(-81/512, 2.641)--(0, 2.366)--cycle;
\filldraw[draw=black, fill=yellow] (81/256, 2.366)--(81/512, 2.641)--(0, 2.366)--cycle;
\filldraw[draw=black, fill=yellow] (81/128, 2.366)--(0.474, 2.641)--(0.316, 2.366)--cycle;
\filldraw[draw=black, fill=yellow] (-243/512, 2.641)--(-0.356, 2.847)--(-243/1024, 2.641)--cycle;
\filldraw[draw=black, fill=yellow] (-243/1024, 2.641)--(-0.119, 2.847)--(0, 2.641)--cycle;
\filldraw[draw=black, fill=yellow] (243/1024, 2.641)--(0.119, 2.847)--(0, 2.641)--cycle;
\filldraw[draw=black, fill=yellow] (243/512, 2.641)--(0.356, 2.847)--(243/1024, 2.641)--cycle;
\node (a9) at (0,3) [acteur][label = below: ]{};
\node (a10) at (0,3.1) [acteur][label = below: ]{};
\node (a11) at (0,3.2) [acteur][label = below: ]{};
\end{tikzpicture}
\caption{The case $m=8$. Each layer is partitioned into seven equal equilateral triangles with the same height as the layer. The picture proves that $\frac{1}{4} + \frac{1}{4}\left(\frac{9}{16}\right) + \frac{1}{4}\left(\frac{9}{16}\right)^2 + \cdots = \frac{4}{7}$.} \label{m=8}
\end{figure}
\end{center}
Unlike Mabry's and Edgar's proofs, Figures \ref{fig:m=3} and \ref{m=8} prove particular cases of Equation \eqref{k1} when $ku\neq s^2$. 

\section{A common idea among the three proofs without words}
At first glance, Figures \ref{fig:Mabry's proof}, \ref{fig: Edgar's}, and \ref{fig: Viewpoints 2000 Group} are little related. However, they are based on the same idea. To illustrate this, we reposition  the Viewpoints 2000 Group's graph as in Figure \ref{fig: reposition}. We shall prove \eqref{k00} for $a = 1$ and $r\in (0,1)$.
This allows us to start with a triangle of base $1$ instead of $a$ as in Figure \ref{fig: Viewpoints 2000 Group}. 

\begin{center}
\begin{figure}[h!]
\begin{tikzpicture}
[acteur/.style={circle, fill=black,thick, inner sep=0.3pt, minimum size=0.01cm}] 
\draw (0,6)--(6,0) node[pos=0, left, xshift=3pt] {\tiny{$\frac{1}{1-\sqrt{r}}$}};
\draw (0,6)--(6,0) node[pos=1, above, xshift=6pt, yshift=-2pt] {\tiny{$\frac{1}{1-\sqrt{r}}$}};
\draw (0,6)--(6,0) node[pos=0, right, yshift=3pt, xshift=-3pt] {\tiny{$\textbf A$}};
\draw (0,6)--(6,0) node[pos=1, below, yshift=2.5pt] {\tiny{$\textbf B$}};
\draw (0,6)--(3.6,0) node[pos=1, below, yshift=2.5pt] {\tiny{$\textbf C$}};
\draw (0,0)--(7,0);
\draw (0,0)--(0,7);
\filldraw[draw=black, fill=yellow] (3.6, 0) -- (6, 0) -- (3.6, 2.4) -- cycle;
\filldraw[draw=black, fill=yellow] (3.6, 2.4) -- (2.16, 3.84) -- (2.16, 2.4) -- cycle;
\filldraw[draw=black, fill=yellow] (1.296,3.84) -- (2.16, 3.84) -- (1.296, 4.704) -- cycle;
\filldraw[draw=black, fill=yellow] (1.296,4.704) -- (0.7776, 4.704) -- (0.7776, 5.2224) -- cycle;
\node (a1) at (0.66,5.1804) [acteur][label = below: ]{};
\node (a2) at (0.54,5.34) [acteur][label = below: ]{};
\node (a3) at (0.42,5.484) [acteur][label = below: ]{};

\draw (3.6,0)--(6,0) node[pos=0.5, below, yshift=2.5pt] {\tiny{$1$}};
\draw (3.6,0)--(3.6,2.4) node[pos=0.5, left, xshift=3pt] {\tiny{$1$}};
\draw (3.6,2.4)--(2.16,2.4) node[pos=0.5, below, yshift=2.5pt] {\tiny{$\sqrt{r}$}};
\draw (2.16, 3.84)--(1.296,3.84) node[pos=0.5, below, yshift=2.5pt] {\tiny{$r$}};
\draw (3.6,2.4)--(6, 0) node[pos=0.5,right] {\tiny{$1$\textsuperscript{st} layer}};
\draw (3.6,2.4)--(2.16,3.84) node[pos=0.5,right] {\tiny{$2$\textsuperscript{nd} layer}};
\draw (1.296,4.704)--(2.16,3.84) node[pos=0.5,right] {\tiny{$3$\textsuperscript{rd} layer}};
\end{tikzpicture}
\caption{Viewpoints 2000 Group's graph repositioned.} \label{fig: reposition} 
\end{figure}
\end{center}
We can now apply Mabry's and Edgar's idea to prove the formula for geometric series.  
The shaded area is clearly $$\frac{1}{2}+\frac{1}{2}r+\frac{1}{2}r^2 + \frac{1}{2}r^3 + \cdots.$$ On the other hand, the shaded area covers $(1+\sqrt{r})^{-1}$ of the area of each layer. Hence, the shaded area covers $(1+\sqrt{r})^{-1}$ of the area of $\triangle{\textbf{ABC}}$. Since $\triangle\textbf{ABC}$ has area $\frac{1}{2(1-\sqrt{r})}$, we have
$$\frac{1}{2}+\frac{1}{2}r+\frac{1}{2}r^2 + \frac{1}{2}r^3 + \cdots\ =\ \frac{1}{1+\sqrt{r}}\cdot \frac{1}{2(1-\sqrt{r})},$$
which is simplified to 
$$1+r+r^2+r^3+\cdots \ =\ \frac{1}{1-r}.$$

\section*{Acknowledgement} 
The authors would like to thank the anonymous referees for various helpful suggestions.

\medskip

\noindent MSC2020: 40A05, 97-01

\end{document}